\documentclass[smallextended]{svjour3}       
\smartqed  
\usepackage{amsfonts}
\usepackage{amsmath}
\usepackage{graphicx}
\usepackage{epsfig}
\usepackage{amssymb}
\usepackage{tikz}
\usetikzlibrary{positioning}
\usetikzlibrary{shapes}

\usepackage{comment}
\usepackage{tkz-graph}
\usepackage{tkz-berge}

\newenvironment{PrfFact}{{\bf Proof.}}{{\hfill{$\blacksquare$\\ }}}

\newlength\flitemwidth

\journalname{Journal of Combinatorial Optimization}

\begin{document}

\title{Graphs with multi-$4$-cycles and the Barnette's conjecture }

\author{Jan~Florek}

\institute{J.~Florek \at  Faculty of Pure and Applied Mathematics,
 Wroclaw University of Science and Technology,
 Wybrze\.{z}e Wyspia\'nskiego 27,
50--370 Wroc{\l}aw, Poland\\
\email{jan.florek@pwr.edu.pl}}

\date{Received: date / Accepted: date}

\maketitle

\begin{abstract}
Let ${\cal H}$ denote the family of all graphs with multi-$4$-cycles and suppose that $G \in {\cal H}$. Then, $G$ is a bipartite graph with a vertex bipartition $\{V_{\alpha}, V_{\beta}\}$. We prove that for every vertex $v \in V_{\beta}$ and for every $2$-colouring $V_{\alpha} \rightarrow \{1, 2\}$ there exists a $2$-colouring $V_{\beta} \rightarrow \{1, 2\}$ such that every cycle in $G$ is not monochromatic and $b(v) = 1$ ($b(v) = 2$).

Let now $G$ be a simple even plane triangulation with a  vertex $3$-partition $\{V_{1}, V_{2}, V_{3}\}$. Denote by $B_{i}$, $i = 1, 2, 3$, the set of all vertices in $V_i$ of degree at least $6$ in $G$. Suppose that $G[B_{1}\cup B_{3}]$ ($G[B_{2}\cup B_{3}]$) is a subgraph of $G$ induced by the set $B_{1}\cup B_{3}$ ($B_{2}\cup B_{3}$, respectively). Let  $G^{*}$ be the dual graph of $G$ with the following $3$-face-colouring: a face $f$ of $G^{*}$ is coloured with $i$ if and only if the vertex $v = f^{*} \in V_{i}$.  We prove that if $H = G[B_{1}\cup B_{3}] \cup G[B_{2}\cup B_{3}] \in {\cal H}$, then, for  any edge chosen on a face  coloured $3$ and of size at least $6$ in $G^{*}$, there exists a Hamilton cycle of $G^{*}$ which avoids this edge.  Moreover, if  every component of $H$ is $2$-connected, then  there exists a Hamilton cycle of $G^{*}$ such that for every face coloured $3$ it avoids every second edge of this face  or it avoids at most two edges of this face.
 \end{abstract}




\section{Introduction}

All graphs considered in this paper will be finite and simple. We use \cite{flobar3} as a reference for undefined terms. In particular, $V(G)$  is the vertex set and $E(G)$ is the edge set of a graph $G$. $N(v)$ is the set of all neighbours of a vertex $v$ in~$G$. The degree $d_{G}(v)$ of a vertex $v$ is the number of edges at $v$. If $\{S, T\}$ is a partition of $V(G)$, then  $E[S, T]$ is the set of all  edges of~$G$ with one end in $S$ and the other end in $T$. If $U \subseteq V(G)$, then $G[U]$ denotes a  subgraph of~$G$ induced by $U$.

In $1884$, P. G. Tait posed a conjecture that all cubic $3$-connected plane graphs had a Hamilton cycle. The conjecture was of special interest since it implied the four colour theorem. Because of counterexamples found the conjecture has been modified in various manners.

Let ${\cal B}$ denote the family of all $3$-connected cubic bipartite plane graphs. In $1969$, Barnette  (\cite{flobar1}) conjectured that every graph in ${\cal B}$ has a Hamilton cycle. Goodey \cite{flobar5} proved that if a graph in ${\cal B}$ has only faces with $4$ or $6$ sides, then it is hamiltonian. 
Bagheri, Feder, Fleischner and Subi \cite{flobar2} study the existence of hamiltonian cycles in plane cubic graphs having a facial $2$-factor. The problem whether a cubic bipartite planar graph has a Hamilton cycle (without the assumption of $3$-connectivity) is NP-complete, as shown by Takanori, Takao and Nobuji \cite{flobar11}. Holton, Manvel and McKay \cite{flobar8} used computer search to confirm the Barnette's conjecture for graphs in ${\cal B}$ with up to $64$ vertices and  they also verifed the following properties $H^{+-}$ and $H^{--}$ for graphs with up to $40$ vertices:
\\[4pt]
$H^{+-}$:  If any two edges are chosen on the same face, then there is a Hamilton cycle through one and avoiding the other,
\\[4pt]
$H^{--}$:  If any two edges are chosen which are an even distance apart on the same face, then there is a~Hamilton cycle that avoids both.
\\[4pt]
Kelmans \cite{flobar9}  proved that Barnette's conjecture holds if and only if every graph in ${\cal B}$ has property $H^{+-}$ (called also Kelmans property).

Stein  expresses hamiltonicity in terms of the dual graph.  Let $G$ be a simple plane triangulation. Stein \cite{flobar10} proved that  there exists a  partition  $\{S, T\}$ of $V(G)$ such that each of partition sets induces an acyclic subgraph in $G$  if and only if $\{ e^{*}\in E(G^{*}): e \in E(S, T)\}$ is the  set of edges of a Hamilton cycle in the dual graph $G^{*}$ (see also Hakimi and Schmeichel~\cite{flobar7})

Let $\cal E$ be the dual family to $\cal B$. Thus, $\cal E$ is the family of all simple even plane triangulations with at least four vertices. Let $G \in \cal E$ and suppose that $\{V_{1}, V_{2}, V_{3}\}$ is a vertex $3$-partition of $V(G)$. We say that a vertex $v$ is \textsl{big} (\textsl{small}) in $G$ if $d_{G}(v) \geq 6$ ($d_{G}(v) = 4$, respectively). Denote by $B_i$ the set of all big  vertices in $V_i$, $i = 1, 2, 3$.
Let ${\cal P}_{G}$ denote the family of all paths $P$ of length at least $2$ in $G$,  with big ends and small inner vertices such that either~$P$ is an induced path in $G$, or the ends $a$, $b$ of $P$ are adjacent and $P + ab$ is an induced cycle in $G$. Florek~\cite{flobar4} proved the following lemma:
\begin{lemma}\label{lemma1.1}
Suppose that $G \in {\cal E}$ has at least three big vertices, and that $X$, $Y$ are disjoint sets of vertices in the graph $G$ satisfying:

\mbox{\rm{(a)}} $B_1 \subseteq X$, $B_2 \subseteq Y$, and $B_3 \subseteq X \cup Y$,

\mbox{\rm{(b)}} for every path in ${\cal P}_G$ the set of its inner vertices is either contained in or is disjoint with $X \cup Y$.

Assume also that the induced graphs $G[X]$ and $G[Y]$ are acyclic.
Then it is possible to partition the vertex set of $G$
into two subsets $S$, $T$ so that each of them induces a tree, and $X \subseteq S$, $Y \subseteq T$.
\end{lemma}
Hence, by the Stein theorem, $\{ e^{*}\in E(G^{*}): e \in E(S, T)\}$ is the  set of edges of a Hamilton cycle in the dual graph $G^{*}$.  Moreover, Florek \cite{flobar4} proved that if subgraphs $G[B_{1}\cup B_{3}]$ and $G[B_{2}\cup B_{3}]$ are acyclic in $G$, then  $G^{*}$ has $H^{+-}$ and $H^{--}$ properties.

Let ${\cal H}$ denote the family of all graphs which contain only cycles of length congruent to $0$ mod $4$. Let $G \in {\cal H}$ and suppose that  $\{V_{\alpha}(G), V_{\beta}(G)\}$ is a vertex bipartition of $V(G)$.  In the first part of our paper we study the structure of graphs in  ${\cal H}$ (see Theorem  \ref{theorem2.1}), and  we prove the following theorem Theorem~\ref{theorem2.2}: For every vertex $v \in V_{\beta}(G)$ and for  every $2$-colouring $a: V_{\alpha}(G) \rightarrow \{1, 2\}$ there  exists a $2$-colouring $b: V_{\beta}(G) \rightarrow \{1, 2\}$ such that every cycle in $G$ is not monochromatic and  $b(v) = 1$ ($b(v) = 2$).

In the second part of our paper we show  applications of  Theorems \ref{theorem2.1}--\ref{theorem2.2} to the Barnette's Conjecture.  Let $G \in \cal E$ and let $H = G[B_{1}\cup B_{3}] \cup G[B_{2}\cup B_{3}]$. First we prove the following Theorem~\ref{theorem2.3}: If  $H \in \cal H$, then for every edge $vw$ with $v \in B_{3}$  and $w \in V_1$$ ($$w \in V_2$$)$ there exists a partition $\{S, T\}$ of $V(G)$  so that each of partition sets induces a tree-subgraph of $G$, $B_{1} \subseteq S$, $B_{2} \subseteq T$, and $v, w \in S$ $($$v, w \in T$, respectively$)$.

Let us define the following $3$-face-colouring of $G^{*}$: a face $f$ of $G^{*}$ is coloured with $i$ if and only if the vertex $v = f^{*} \in V_{i}$. Then, we obtain a dual version of Theorem \ref{theorem2.3}: If $H \in \cal H$, then, for  any edge chosen on a face  coloured $3$ and of size at least $6$ in $G^{*}$, there exists a Hamilton cycle of~$G^{*}$ which avoids this edge.

Finally, we prove the following Theorem~\ref{theorem2.4}: If $H \in \cal H$ and if every component of $H$ is $2$-connected, then there exists a partition $(S, T)$ of $V(G)$ so that each of partition sets induces a tree-subgraph in $G$, $B_{1} \subseteq S$, $B_{2} \subseteq T$, and for every $v \in B_{3}$  the following implications are satisfied:
\begin{enumerate}
\item [($1$)] if $d_{H}(v) \geqslant 3$, then  $N(v) \cap V_1 \subseteq S$  and $N(v) \cap V_2 \subseteq T$,
\item [($2$)]
if  $d_{H}(v) = 2$ and $v \in S$ ($v \in T$), then at most two vertices of  $N(v)\cap V_1$ ($N(v)\cap V_2$) belong to  $S$ (or $T$)  and all other neighbours of $v$ belong to $T$ ($S$, respectively).
\end{enumerate}
Then, we obtain a dual version of Theorem \ref{theorem2.4}: If $H \in {\cal H}$ and if  every component of $H$ is $2$-connected, then  there exists a Hamilton cycle of $G^{*}$ such that for every face coloured $3$ it avoids every second edge of this face  or it avoids at most two edges of this face.

\section{Graphs with multi-$4$-cycles}

Let ${\cal H}$ denote the family of all graphs such that each cycle has a length congruent to $0$ mod $4$. Suppose that $G \in \cal H$ and assume that $\{V_{\alpha}(G), V_{\beta}(G)\}$ is a vertex bipartition of $G$.  Vertices belonging to $V_{\alpha}$ (or $V_{\beta}$) are called of \textit{type}~$\alpha$ (\textit{type} $\beta$ respectively). For any $2$-colouring $a: V_{\alpha}(G) \rightarrow \{1, 2\}$ and for any $2$-colouring $b: V_{\beta}(G) \rightarrow \{1, 2\}$ we define a $2$-colouring $a \triangledown b:  V_{\beta}(G) \rightarrow \{1, 2\}$ (not necessarily proper) such that $a \triangledown b(w) = a(w)$, for $w \in V_{\alpha}(G)$, and $a \triangledown b(w) = b(w)$, for $w \in V_{\beta}(G)$.
A cycle in $G$ is \textit{monochromatic} with respect to $a \triangledown b$ if all vertices of this cycle have the same colour. Given a subgraph $C$ of $G$, we call a path $P$ a $C$-\textsl{path} if $P$ is non-trivial and meets $C$ exactly in its ends.
Let $P = x_{0} x_{1} \ldots  x_{k}$. For $0 \leqslant i \leqslant j \leqslant k$ we write  (see Diestel \cite{flobar3})
$$Px_{i} =  x_{0}\ldots  x_{i}$$
$$x_{i}P =  x_{i} \ldots  x_{k}$$
$$x_{i}Px_{j} = x_{i} \ldots  x_{j}$$
for the appropriate subpaths of $P$. We use a similar notation for the concatenation of paths; for example, if the union  $xPy \cup yQx$ of two paths is a cycle, we may simply denote it by $xPyQx$.

\begin{definition}\label{def2.1} If $P$ is a path of length at least $2$, then $Int P$ is the set of its all inner vertices. If $P$ has length $1$, then $Int P$ is the edge of $P$.
\end{definition}
\begin{definition}\label{def2.2}
Let $G \in {\cal H}$ and suppose that $\{V_{\alpha}(G), V_{\beta}(G)\}$ is a vertex bipartition. We say that a pair $(P, Q)$ of two disjoint  paths in $G$  \textsl{cuts} $G$ if the following conditions are satisfied:

(a) $P$ (and $Q$) has all inner vertices of degree $2$ in $G$, and has ends  of different types which are of degree at least $3$ in $G$,

(b) the graph $G - (Int P \cup Int Q)$ has two components, which are called subgraphs of $G$ \textsl{determined by $(P, Q)$}, and each of them contains ends of these paths of  the same type.
\end{definition}
\begin{lemma}\label{lem2.1}
Let $G \in {\cal H}$ and suppose that $\{V_{\alpha}(G), V_{\beta}(G)\}$ is a vertex bipartition. If $C$ is a $2$-connected subgraph of $G$ and $P$ is a $C$-path, then ends of $P$ are of the same type.
\end{lemma}
\begin{PrfFact}
Let $C$ be a $2$-connected subgraph of $G$ and suppose that $Q$, $R$ are independent $x$-$y$ paths  contained in $C$. Let  $P$ be a $C$-path connecting $x$ and~$y$. Notice that if $x$ and $y$ are of different types, then one of cycles $xPyQx$, $xPyRx$ or $xQyRx$ has odd number of vertices of type $\alpha$ (and also of  type $\beta$). Hence, it has $4n +2$ vertices, which is a contradiction.
\end{PrfFact}
\begin{lemma}\label{lem2.2}
Let $G \in {\cal H}$ and suppose that $\{V_{\alpha}(G), V_{\beta}(G)\}$ is a vertex bipartition. If $G$ is $2$-connected, then for every $x$-$y$ path $P$ in $G$ satisfying condition (a) of Definition \ref{def2.2} there exists a sequence  $B_{1}, B_{2},\ldots,  B_{n+1}$, $n \geqslant 1$, of blocks in $G-Int P$ and a sequence $a_{1},a_{2},\ldots,  a_{n}$ of cut vertices, such that $G-Int P = B_{1} \cup \ldots  \cup B_{n+1}$, blocks $B_1$ and $B_{n+1}$ are $2$-connected, $x \in B_1$, $y \in B_{n+1}$,  and $\{a_{i}\} = B_{i} \cap B_{i+1}$, $a_{1} \neq x$,  $a_{n} \neq y$.
\end{lemma}
\begin{PrfFact}
Since $G$ is $2$-connected and  $P$ has ends of different types,  by Lemma~\ref{lem2.1}, $G-Int P$ is not $2$-connected. Since  $G$ is $2$-connected and inner vertices of~$P$ are of degree $2$ in $G$, there is a path $Q$ contained in  $G-Int P$ which is connecting $x$ and $y$. Since $G-Int P$ is not $2$-connected, there exists  a sequence $B_{1}, B_{2},\ldots,  B_{n+1}$, $n \geqslant 1$, of  blocks in $G-Int P$ such that $x \in B_1$, $y \in B_{n+1}$, every block of this sequence contains an edge of $Q$  and every edge of $Q$ belongs to a one of them. Thus, there exists a sequence $a_{1},a_{2},\ldots,  a_{n}$ of cut vertices such that $\{a_{i}\} = B_{i} \cap B_{i+1}$.

Since $G$ is $2$-connected every vertex  of $G-Int P$ belongs to a path contained in  $G-Int P$ which is connecting $x$ and $y$, whence it is connecting two vertices of a block $B_{i}$, for some $i = 1, \ldots, n+1$. Thus, $G-Int P = B_{1} \cup \ldots  \cup B_{n+1}$.

Certainly, $x \neq a_1$ and $y \neq a_n$, because $G$ is $2$-connected. Notice that blocks $B_{1}$ and $B_{n+1}$ are $2$-connected, because $x$, $y$ are of degree at least $3$ in $G$.
\end{PrfFact}
\begin{lemma}\label{lem2.3}
Let $G \in {\cal H}$ and suppose that $\{V_{\alpha}(G), V_{\beta}(G)\}$ is a vertex bipartition.  If $G$ is $2$-connected, then no $4$-cycle  has two adjacent vertices of degree at least~$3$ in $G$.
\end{lemma}
\begin{PrfFact}
Let $abcda$ be a $4$-cycle in $G$ and suppose that ends of the edge $e = ad$ are of degree at least $3$ in $G$.  Since ends of $e$ are of different types, by Lemma~\ref{lem2.2}, $G-e = B_{1} \cup B_{2}$ and one of the following cases occurs:

(i) $B_1$ is $2$-connected,  $B_2$ is connected, $B_{1}\cap B_2 = b$, $ab$ is an edge of $B_1$ and $bcd$ is a path in $B_{2}$,

(ii) $B_1$ is connected, $B_2$ is $2$-connected $B_{1}\cap B_{2} = c$, $abc$ is a path in $B_{1}$ and $cd$ is an edge of $B_2$.

In the case (i),  $adcb$ is a $B_{1}$-path. Hence, by Lemma \ref{lem2.1}, vertices $a, b$ are of the same type, but it is impossible, because $ab$ is an edge. Similarly, we obtain a contradiction in the case (ii). Hence, $e$ has an end-vertex of degree $2$.
\end{PrfFact}
\begin{lemma}\label{lem2.4}
Let $G \in {\cal H}$ be $2$-connected  and suppose that $\{V_{\alpha}(G), V_{\beta}(G)\}$ is a vertex bipartition. For every path $P$ in $G$ satisfying condition (a) of Definition~\ref{def2.2} and for every $2$-connected block $B$ in $G - Int P$, there exists a path $Q$ in $G$ such that $(P, Q)$ cuts $G$ and $B$ is contained in a component of $G - (Int P \cup Int Q)$.
\end{lemma}
\begin{PrfFact}
Let  $P$ be a path of $G$ satisfying condition (a) of Definition \ref{def2.2} and suppose that $x$ (or $y$) is its end of type $\alpha$ (or $\beta$, respectively).
By Lemma~\ref{lem2.2}, there exists a sequence  $B_{1}, B_{2},\ldots,  B_{n+1}$, $n \geqslant 1$, of blocks in $G-Int P$ and a sequence $a_{1},a_{2},\ldots,  a_{n}$ of cut vertices, such that $G-Int P = B_{1} \cup \ldots  \cup B_{n+1}$,  $B_1$ and $B_{n+1}$ are $2$-connected, $\{a_{i}\} = B_{i} \cap B_{i+1}$ and $x \in B_1$, $x \neq a_1$ and $y \in B_{n+1}$, $y \neq a_n$. Since $B_1$ and $B_{n+1}$ are $2$-connected, by  Lemma \ref{lem2.1}, vertices $a_1$ and $x$ ($a_{n}$ and $y$) are of type $\alpha$ ($\beta$, respectively). We can assume that $B = B_l$. Since $B_l$ is $2$-connected, by  Lemma  \ref{lem2.1}, vertices $a_{l-1}$, $a_l$ are of the same type. We assume that these vertices are of type $\beta$ (if they are of type $\alpha$, the proof is analogous).

 Let $1 < k \leqslant n$ be the first integer such that $d_{G}(a_{k}) \geqslant 3$ and $a_k$ is of type~$\beta$, and suppose that $1 \leqslant  j < k$ is the last integer such that $d_{G}(a_{j}) \geqslant 3$.  Then, ${Q =  a_{j}a_{j+1}\ldots a_{k}}$ is a path satisfying condition (a) of Definition \ref{def2.2}.
Notice that paths $P$, $Q$ are disjoint and $G - (Int P \cup Int Q)$ has two components $C = B_{1} \cup \ldots  \cup B_{j}$ and $D = B_{k+1} \cup \ldots  \cup B_{n+1}$. Since $x, a_j \in C$ are of type $\alpha$ and $a_k, y \in D$ are of type $\beta$, the pair $(P, Q)$ cuts $G$. It follows, by definition of $k$, that $k \leqslant l-1$, because $d_{G}(a_{l-1}) \geqslant 3$ and  $a_{l-1}$ is of type $\beta$. Hence, $B = B_l \subseteq D$.
\end{PrfFact}
\begin{theorem}\label{theorem2.1}
Let $G \in {\cal H}$ be $2$-connected and suppose that  $\{V_{\alpha}(G), V_{\beta}(G)\}$ is a vertex bipartition.  If $C$  is a minimal subgraph in $G$ determined by a pair of paths which cuts $G$, then its any two vertices  of degree at least $3$  in $G$ are of the same type.
\end{theorem}
\begin{PrfFact}
Let $C$  be a minimal subgraph of $G$ determined by a pair of paths (say $(P, Q)$) which cuts $G$. Hence, $G - (Int P \cup Int Q)$ consists of two components $C$ and $D$. Without loss of generality we can assume that ends of paths $P$, $Q$ belonging to $C$ are of type $\alpha$.  If $C$ is not $2$-connected, then, by Lemma  \ref{lem2.2}, there exists a sequence  $B_{1}, B_{2},\ldots,  B_{n+1}$ of blocks  in $C$, $n \geqslant 1$ and a sequence $a_{1},a_{2},\ldots,  a_{n}$ of cut vertices  such that $C = B_{1} \cup \ldots  \cup B_{n+1}$ and $\{a_{i}\} = B_{i} \cap B_{i+1}$ (if $C$ is $2$-connected, then we put $C = B_1$).

Assume that $C$ has two vertices, say $x$ and $y$, of different types which are of degree at least $3$ in $G$. We shall prove that $C$ is not minimal. Without loss of generality we can assume that one of the following conditions occurs:

(i) $x =  a_{i}$, $y =a_{k}$, for some $i < k$, and  $B_j$ is an edge, for every $i < j \leqslant k$,

(ii)  $x$ and $y$ belong to the same $2$-connected block, say $B_l$.

Case (i). Then, there exists a path  $R = a_{i}a_{i+1}\ldots a_{k}$ satisfying condition (a) of Definition \ref{def2.2}. If $a_i$ is of type $\alpha$ and $a_{k}$ is of type $\beta$, then $(P, R)$ cut~$G$, $B_{1} \cup \ldots \cup  B_{i}$ and $B_{k+1} \cup \ldots \cup B_{n+1} \cup Q \cup D$ are two components of $G - (Int P \cup Int R)$.  Similarly,  if $a_i$ is of type $\beta$ and $a_{k}$ is of type $\alpha$, then $(R, Q)$ cut $G$, $B_{k+1} \cup \ldots \cup B_{n+1}$ and $D \cup P \cup B_{1} \cup \ldots \cup  B_{i}$ are two components of $G - (Int R \cup Int Q)$. Hence, $C$ is not minimal.

Case (ii). Without loss of generality we can assume that vertices $a_{l-1}$ and $a_{l}$ are of type $\alpha$, because, by Lemma \ref{lem2.1}, they are of the same type. Since $x$ and $y$ are of different type, there exists a path $P_{1} \subset B_l$ satisfying condition (a) of Definition \ref{def2.2} and one of its ends, say $x_1$, is of type $\beta$. Hence, there exists a path $S$ in $B_l$ connecting $a_{l-1}$ and $a_{l}$ and omitting the vertex $x_1$, because $B_l$ is $2$-connected. Notice that $S$ is disjoint with $Int P_1$. Since $G$ is $2$-connected, there exists a path $T$  in $G$ connecting  $a_{l-1}$ and $a_{l}$  which  contains paths $P$ and $Q$. Let~$B$ be a $2$-connected block of $G-Int P_{1}$ which contains a cycle $a_{l-1}Sa_{l}Ta_{l-1}$. Since every vertex of $G-B_{l}$ belongs to a path in $G -Int P_{1}$ connecting $a_{l-1}$ and $a_{l}$, $G-B_{l} \subset B$. Hence, $G - B \subset B_l$. By Lemma~\ref{lem2.4}, there exist a path $Q_1$ of $G$  such that $(P_{1}, Q_{1})$ cuts $G$ and one component of $G-(Int P_{1} \cup Int Q_{1})$, say $E$, is disjoint with $B$. Thus, $E \subseteq G - B \subset B_l \subseteq C$.  Hence, $C$ is not minimal.
\end{PrfFact}
\begin{theorem}\label{theorem2.2}
Let $G \in {\cal H}$ and suppose that  $\{V_{\alpha}(G), V_{\beta}(G)\}$ is a vertex bipartition.  For every vertex $v \in V_{\beta}(G)$ and for  every $2$-colouring $a: V_{\alpha}(G) \rightarrow \{1, 2\}$ there  exists a $2$-colouring $b: V_{\beta}(G) \rightarrow \{1, 2\}$ satisfying the following conditions:

(1)  every cycle in $G$ is not monochromatic with respect to $a \triangledown b$,

(2)  if a path in $G$ has ends of type $\beta$ and its all inner vertices are of degree~$2$ in $G$, then its vertices of type $\beta$ are coloured $1$ and $2$ alternately by $b$,

(3) $b(v) = 1$ $($$b(v) = 2$$)$.
\end{theorem}
 \begin{PrfFact}
Let $G \in {\cal H}$. Without loss of generality we can assume that $G$ is connected. We apply induction on $|G|$. Fix a vertex $v \in V_{\beta}(G)$ and fix a $2$-colouring $a: V_{\alpha}(G)\rightarrow \{1, 2\}$.

Assume first that $G$ is not $2$-connected. Pick a block $B_0$ of $G$ such that $v \in V_{\beta}(B_{0})$.  Let $B_{1}, \ldots, B_{n}$ be a sequence of  connected subgraphs of $G$ and let $v_1,\ldots, v_{n}$ be a sequence of vertices such that $v_{i}$ is the unique
 vertex of $B_{0} \cap B_{i}$, and subgraphs $B_{i} \backslash \{v_{i}\}$ are pairwise disjoint.
By induction (we apply Theorem~\ref{theorem2.2} for the graph $B_{i}$ and  the vertex $v_{i}$, $i = 1,\ldots, n$, and  for the graph $B_{0}$ and the vertex $v$), there exists a  $2$-colouring $b_{i}: V_{\beta}(B_{i})\rightarrow \{1, 2\}$ satisfying the assertions:
\begin{enumerate}
\item [($\circ$)]
every cycle in $B_i$ is not monochromatic in $a \triangledown b_{i}$,
\item[($\circ$)]
 if a path has ends of type $\beta$ and has all  inner vertices of degree $2$ in $B_i$,
then its all vertices of type $\beta$ are coloured $1$ and $2$ alternately by $b_i$,
\item [($\circ$)]
if $v_{i}$ is of type $\beta$, then $b_{i}(v_{i}) = b_{0}(v_{i})$,  for $i \geqslant 1$,
\item [($\circ$)]
 $b_{0}(v) = 1$ $($$b_{0}(v) = 2$$)$.
\end{enumerate}
Notice that a $2$-colouring $b = b_{0}\triangledown \ldots \triangledown b_{n}$ satisfies assertions (1)-(3) of Theorem~\ref{theorem2.2}.  Thus, it is sufficient to prove Theorem \ref{theorem2.2} for $2$-connected graphs in~${\cal H}$.

First we prove the following:
\begin{enumerate}
\item [($A$)] If $K$ is  a connected subgraph of $G$ such that its all vertices of degree at least $3$  in $K$ are of  type $\beta$, then there exists a $2$-colouring $k : V_{\beta}(K) \rightarrow \{1, 2\}$ so that
\item [($a$)]
every cycle in $K$ is not monochromatic with respect to $a \triangledown k$,
\item [($b$)]
if a path in $K$ has ends of type $\beta$, then its vertices of type $\beta$ are coloured $1$ and $2$ alternately by $k$,
\item [($c$)]
if $v \in V_{\beta}(K)$, then $k(v) = 1$  $($$k(v) = 2$$)$.
\end{enumerate}
Proof of ($A$). Fix a vertex $w \in V_{\beta}(K)$. For every $u \in  V_{\beta}(K)$ there exists a path  $P$ in $K$ connecting $w$ and $u$. Ler $k(w) = 1$ and assume that vertices of $V_{\beta}(P)$ are coloured, by $k$, with $1$ and with $2$ alternately. We prove that the colouring of $u$ is  independent of the choice of the path in $K$ connecting $w$ and $u$.

Let $P$ and $Q$ be any two paths in $K$ connecting $w$ and $u$. Let $d(P)$ (or $d(Q)$) denote length of the path $P$ (or $Q$, respectively).  Notice that $P$ (or $Q$) has ends of the same colour if and only if $d(P)$ ($d(Q)$, respectively) is congruent to $0$ modulo $4$.  If $P$ and $Q$ are independent, then $d(P)$ and $d(Q)$ are congruent modulo  $4$, because $G \in {\cal H}$. Hence follows that the colouring of $u$ is  independent of the choice of the path $P$ or $Q$.

If $P$ and $Q$ are not independent, then their common inner vertices are of degree at least $3$ in $K$. Hence, they are of type $\beta$. Then there are two sequences of paths $P_{1},\ldots,  P_{n}$ and $Q_{1},\ldots,  Q_{n}$  such that $P = P_{1} \cup \ldots \cup  P_{n}$, $Q = Q_{1} \cup \ldots \cup  Q_{n}$ and,  for every $i = 1,\ldots n$,  the paths $P_i$, $Q_i$ are independent and have common ends of type $\beta$, or  $P_i =Q_i$. Thus, the colouring $k$ is compatible for all vertices of $V_{\beta}(K)$. Certainly, if $v \in V_{\beta}(K)$, then we may interchange the colours $1$ and $2$ in $V_{\beta}(K)$ so that condition ($c$) is satisfied. Thus, conditions ($b$)-($c$) hold. Certainly, ($a$) follows from ($b$), which completes the proof of ($A$).

Now we prove the following:
\begin{enumerate}
\item[($B$)]  If $L$ is a connected subgraph of $G$ such that its all vertices of degree at least $3$  in $G$ are of type $\alpha$, then there exists a $2$-colouring $l : V_{\beta}(L) \rightarrow \{1, 2\}$ so that
\item [($d$)]
every cycle in $L$ is not monochromatic with respect to $a \triangledown l$,
\item [($e$)]
 if a path in $L$ has ends of type $\beta$ and its all  inner vertices are of degree $2$ in $G$, then its all vertices of type $\beta$ are coloured $1$ and $2$ alternately by $l$,
\item [($f$)]
every path in $L$ not containing $v$ and with ends of degree at least $3$  in $G$ is not monochromatic with respect to $a \triangledown l$,
\item [($g$)]
if  $v \in V_{\beta}(L)$, then $l(v) = 1$  $($$l(v) = 2$$)$.
\end{enumerate}
Proof of ($B$). Let ${\cal L}_1$ be the set of all paths in $L$ with ends of type $\beta$ and with all  inner vertices of degree $2$ in $G$. If $P \in {\cal L}_1$, then  we assume that vertices of $V_{\beta}(P)$ are coloured, by $l$, with $1$ and  with $2$ alternately.  Certainly, if $v \in V_{\beta}(P)$, then we can interchange colours $1$ and $2$ in $V_{\beta}(P)$ so that condition ($g$) is satisfied. Notice that that the colouring $l$ is compatible for all vertices of the union of $V_{\beta}(P)$, for $P \in {\cal L}_1$. Hence, conditions  ($e$) and ($g$) hold.

Let ${\cal L}_2$ be the set of all paths in $L$ of length $2$ with ends of degree at least~$3$ in $G$ and not containing $v$. If $P \in {\cal L}_2$, then we assume that the inner vertex of $P$ (of type $\beta$) is coloured, by $l$, so that $P$ is not monochromatic. If $P$ is a path in $L$ of length at least $4$ and with ends of degree at least $3$ in $G$ (of type $\alpha$), then it contains a path belonging to ${\cal L}_{1}$ or ${\cal L}_{2}$. Hence, it is not monochromatic with respect to $a \triangledown l$. Thus, condition  ($f$) holds.

 Certainly, condition ($d$) follows from conditions ($e$)-($f$) which completes the proof of ($B$).

Now we assume that $G$ is $2$-connected. If all vertices of degree at least~$3$ in $G$ are of  the same type, then, by ($A$) (or ($B$)), there exists a $2$-colouring $k : V_{\beta}(G) \rightarrow \{1, 2\}$ ($l : V_{\beta}(G) \rightarrow \{1, 2\}$) satisfying  conditions  ($a$)-($c$) (($d$)-($f$), respectively). Hence, Theorem \ref{theorem2.2} holds.

Thus, we can assume that $G$ has a path which satisfies condition (a) of Definition \ref{def2.2}. By Lemma  \ref{lem2.4} and Theorem \ref{theorem2.1}, there exists a pair $(P, Q)$ of paths which cuts $G$, $G-(Int P \cup Int Q)$ is the union of two componets, say $C$ and $D$, and all  vertices  of $C$ which have degree at least $3$  in $G$ are of the same type. Let us denote by $x_{1},y_{1}$ ($x_{2},y_{2}$) ends of $P$ ($Q$, respectively).  We may assume that $x_{1}, x_{2} \in C$. Let us consider the following cases:

(i)  all  vertices  of $C$ which have degree at least $3$  in $G$ are of type $\beta$,

(ii) all  vertices  of $C$ which have degree at least $3$  in $G$ are of type $\alpha$.

Case (i). By induction (for the graph $D$), there exists a $2$-colouring $d : V_{\beta}(D) \rightarrow \{1, 2\}$ satisfying assertions  (1)-(2) and (3) (if $v \in  V_{\beta}(D)$).  

Let $K = C \cup P \cup Q$. Notice that all vertices of degree at least $3$ in $K$ are of type $\beta$.  Then, by ($A$)  there exists a $2$-colouring $k : V_{\beta}(K) \rightarrow \{1, 2\}$ satisfying conditions  ($a$)-($b$) and ($c$) (if $v \in  V_{\beta}(K) $). Since vertices $x_{1}$ and $x_{2}$ are of type $\beta$, by condition ($b$), every path in $K$ with ends $x_{1}$ and $x_{2}$ is not monochromatic with respect to $a\triangledown k$. Hence, every cycle containing $P\cup Q$ is not monochromatic with respect to $a\triangledown k\triangledown d$.  Thus, in this case, Theorem~\ref{theorem2.2} holds.

Case (ii).
By ($B$) (for $L = C$), there exists a $2$-colouring  $l: V_{\beta}(C) \rightarrow \{1, 2\}$ satisfying conditions ($d$)-($f$)  and ($g$) (if $v \in V_{\beta}(C)$).

Suppose first that $v\in D \cup P \cup Q$. By induction (for the graph $D \cup P \cup Q$ and for $v$), there exists a $2$-colouring $c_{1} : V_{\beta}(D \cup P \cup Q) \rightarrow \{1, 2\}$ satisfying  assertions  (1)-(3) of Theorem \ref{theorem2.2}. Since $v \notin C$ and vertices $x_1, x_{2}$  have degree at least $3$ in $G$, by condition  ($f$), every path with ends $x_{1}$ and $x_{2}$ contained in $C$ is not monochromatic with respect to $a \triangledown l$. Hence, every cycle containing $P\cup Q$ is not monochromatic with respect to $a \triangledown l \triangledown c_{1}$. Thus, in this case, Theorem~\ref{theorem2.2} holds.

Let now $v \in C$. By induction (for $D \cup P \cup Q$ and for $y_{1} \in V_{\beta}(D \cup P \cup Q)$), there exists a $2$-colouring $c_{2} : V_{\beta}(D \cup P \cup Q) \rightarrow \{1, 2\}$ satisfying assertions  (1)-(2) of Theorem \ref{theorem2.2} and $c_{2}(y_{1}) \neq a(x_{1})$. Hence, every cycle containing $P\cup Q$  is not monochromatic with respect to $a \triangledown l \triangledown c_{2}$. Thus, in this case, Theorem~\ref{theorem2.2} holds.
 \end{PrfFact}
\begin{example}\label{example1.1}
Let $K_{3, 4} = E^{3}\ast E^4$ be a complete bipartite graph, where $E^3$ and $E^4$ are disjoint empty graphs of order $3$ and $4$, respectively. Let $a: E^{4} \rightarrow \{1, 2\}$ be a $2$-colouring such that two vertices are coloured $1$ and the other are coloured~$2$. Notice that $K_{3, 4} \notin {\cal H}$ and for every $2$-colouring $b: E^{3} \rightarrow \{1, 2\}$ there exists a cycle monochromatic with respect to $a \triangledown b$.
\end{example}
\begin{remark}\label{remark2.1}
For every  integer $n > 0$, there exists a plane graph $G \in {\cal H}$ with a vertex bipartition  $\{V_{\alpha}(G), V_{\beta}(G)\}$ and a $2$-colouring $a: V_{\alpha}(G) \rightarrow \{1, 2\}$ such that for every $2$-colouring $b: V_{\beta}(G) \rightarrow \{1, 2\}$ there exists a path of length at least $n$ which is monochromatic with respect to $a \triangledown b$.
\end{remark}
\begin{corollary}\label{coro2.1}
Let $G \in {\cal H}$ and suppose that  $\{V_{\alpha}(G), V_{\beta}(G)\}$ is a vertex bipartition.  Let $v, y\in V_{\beta} (G)$ be vertices of a $4$-cycle. For every $2$-colouring $a: V_{\alpha}(G) \rightarrow \{1, 2\}$,  there exists a  $2$-colouring $b: V_{\beta}(G) \rightarrow \{1, 2\}$ which satisfies the following conditions:

(1)  every cycle is not monochromatic  with respect to $a \triangledown b$,

(2)  $b(v) = 1$ and $b(y) = 2$ ($b(v) = 2$ and $b(y) = 1$).
\end{corollary}
\begin{PrfFact}
Let $G \in {\cal H}$. Without loss of generality we can assume that $G$ is connected. Assume that $v, y\in V_{\beta} (G)$ are vertices of a $4$-cycle (say $vxyzv$). Notice that $x, z \in V_{\alpha}(G)$. Pick a $2$-connected block $B$ of $G$ containing $vxyzv$.

Let $a: V_{\alpha}(G) \rightarrow \{1, 2\}$ be any $2$-colouring.  By  Theorem \ref{theorem2.2} there are $2$-colourings $b_{1}: V_{\beta}(B) \rightarrow \{1, 2\}$ and $b_{2}: V_{\beta}(B) \rightarrow \{1, 2\}$) satisfying the following conditions:
\begin{enumerate}
\item [($a$)]
every cycle in $B$ is not monochromatic with respect to $a \triangledown b_1$ and $a \triangledown b_2$,
\item [($b$)]
 if a path of length $2$ in $B$ has ends of type $\beta$ of degree at least $3$ in $B$ and has the  inner vertex of degree $2$ in $B$, then its ends are coloured  with different colours by $b_1$, and also by $b_2$,
\item [($c$)]
$b_{1}(v) = a(x) = b_ {2}(y)$.
\end{enumerate}
We shall prove that there  exists a $2$-colouring $b_{0}: V_{\beta}(B) \rightarrow \{1, 2\}$ satisfying the following conditions:
\begin{enumerate}
\item [($d$)]
every cycle in $B$ is not monochromatic with respect to $a \triangledown b_{0}$,
\item [($e$)]
$b_{0}(y) \neq  b_{0}(v) = a(x)$ ($b_{0}(v) \neq  b_{0}(y) = a(x)$).
\end{enumerate}
Since $B$ is $2$-connected,  by Lemma \ref{lem2.3}, one of the following cases occurs:
\begin{enumerate}
\item[($i$)]
 $d_{B}(v),  d_{B}(y) \geqslant 3$  and $d_{B}(x) = d_{B}(z) = 2$,
\item [($ii$)]
 $d_{B}(v) =  d_{B}(y) = 2$.
\end{enumerate}

Case ($i$) . By conditions ($b$)-($c$),  $b_{1}(y) \neq  b_{1}(v) = a(x)$ and $b_{2}(v) \neq  b_{2}(y) = a(x)$. Hence, $b_{0} = b_{1}$ ($b_{0} = b_{2}$, respectively)  satisfies conditions ($d$)-($e$).

Case ($ii$). If  $a(x) = a(z)$, then, by conditions ($c$) and ($a$),  $b_{1}(v) = a(x) \neq b_{1}(y)$ and $b_{2}(y)= a(x) \neq b_{2}(v)$. Hence,  $b_{0} = b_{1}$ ($b_{0} = b_{2}$, respectively)  satisfies conditions ($d$)-($e$).

If  $a(x) \neq a(z)$, then, by condition ($c$),  $b_{1}(v) = a(x)$ and $b_{2}(y) = a(x)$. Hence, we can assume that  $b_{1}(y) = a(z)$ and $b_{2}(v) = a(z)$.  Otherwise, we can recolour $y$ and $v$, because  $d_{B}(y) =  d_{B}(v) = 2$.  Hence,  $b_{0} = b_{1}$ ($b_{0} = b_{2}$, respectively)  satisfies conditions ($d$)-($e$).

Finally, we can extend the colouring $b_0$ to a $2$-colouring $b$ of  $V_{\beta}(G)$  satisfying conditions ($1$) and ($2$), by the way presented in the proof of Theorem \ref{theorem2.2}.
\end{PrfFact}
\section{Application to the Barnette's Conjecture}
Let $\cal E$ be the family of all simple even plane triangulations. Let $G \in {\cal E}$ and suppose that  ${V_1, V_2, V_3}$ is a vertex $3$-partition of $V(G)$.
If $v$ is a vertex in $G$, then $N_{i}(v)$ is the set of all its neighbours belonging  to $V_i$, for $i = 1, 2, 3$.

 We say that a vertex $v$ is \textsl{big} (\textsl{small}) in $G$ if $d_{G}(v) \geq 6$ ($d_{G}(v) = 4$, respectively). Denote by $B_i$ (or $S_i$) the set of all big (small, respectively) vertices in $V_i$, $i = 1, 2, 3$. Note that $G[B_{1} \cup B_{2} \cup B_{3}]$ is $2$-connected, because its every face is triangle or square.

Let ${\cal P}_{G}$ denote the family of all paths $P$ of length at least $2$ in $G$,  with big ends and small inner vertices such that either $P$ is an induced path in $G$, or the ends $a$, $b$ of $P$ are adjacent and $P + ab$ is an induced cycle in $G$.  Let $V^{0}(P)$ denote the set of two vertices in $G$ so that each of them is adjacent to every vertex of $P$. Let $V^{1}(P)$ denote the set of two ends of $P$.

Let $C^l$ denote a cycle of length $l$, $E^2$ denotes an empty graph of order~$2$. Then the join ${C^4} \ast {E^2}$ is the octahedron. Notice that if $G \ncong {C^{2l}} \ast {E^2}$, for $l \geqslant 2$, then every small vertex in $G$ belongs to a path of  ${\cal P}_{G}$.

Recall that ${\cal H}$ denote the family of all graphs which contain only cycles of length congruent to $0$ mod $4$.
\begin{theorem}\label{theorem2.3}
Let $G \in {\cal {E}}$ and suppose that $G[B_{1} \cup B_{3}] \cup G[B_{2} \cup B_{3}]  \in  {\cal H}$. For every edge $vw$ with $v \in B_{3}$  and $w \in V_1$$ ($$w \in V_2$$)$ it is possible to partition the vertex set of $G$ into two disjoint subsets $S$, $T$ so that each of them induces a tree, $B_{1} \subseteq S$, $B_{2} \subseteq T$, and $v, w \in S$ $($$v, w \in T$, respectively$)$.
\end{theorem}
 \begin{PrfFact}
Let $G \in {\cal {E}}$. If $G \cong {C^{2l}} \ast {E^2}$, for some $l \geqslant 2$, then Theorem \ref{theorem2.3} holds. Therefore, we assume that $G \ncong {C^{2l}} \ast {E^2}$. Hence, every small vertex in $G$ belongs to a path of ${\cal P}_{G}$.

 Let $H = G[B_{1} \cup B_{3}] \cup G[B_{2} \cup B_{3}]  \in  {\cal H}$.
Notice that  $H$ is a bipartite graph with bipartition $(B_1 \cup B_2, B_3)$.  Let $a: B_1 \cup B_2 \rightarrow \{1, 2\}$ be a $2$-~colouring such that vertices of $B_1$ (or $B_2$) are coloured $1$ ($2$, respectively).  Fix an edge $vw$ such that $v \in B_{3}$. By Theorem \ref{theorem2.2}, there exists a $2$-colouring $b: B_{3} \rightarrow \{1, 2\}$  such that every cycle in $G[B_{1} \cup B_{2} \cup B_{3}]$ is not monochromatic with respect to $a \triangledown b$, and $b(v)$ is coloured $1$ (or $2$), for $w \in V_1$ ($w \in V_2$, respectively).

Let $w \in B_1 \cup B_2$ and suppose that
$$X = B_{1} \cup \{u \in B_{3}: b(u) = 1\} \hbox{ and } Y = B_{2} \cup \{u \in B_{3}: b(u) =2\}.$$
Remark, that $X$ an $Y$ are disjoint and the induced graphs $G[X]$ and $G[Y]$ are acyclic. Hence, by Lemma  \ref{lemma1.1},  it is possible to partition the vertex set of $G$ into two subsets $S$, $T$ so that each of them induces a tree, $X \subseteq S$ and $Y \subseteq T$. Certainly, if $w \in B_1$ ($w \in B_2$), then $v, w \in X$ ($v, w \in Y$, respectively).

Assume now that $w \in S_1 \cup S_2$.  Then, $w$ belongs to a path of ${\cal P}_{G}$, say $P_{w}$, such that all vertices of $V^{0}(P_{w})$ are big in $G$, or $V^{0}(P_{w}) = \{v, y\}$, where $y \in S_{3}$. Let $vxyzv$ be a $4$-cycle  consisting of vertices belonging to $V^{0}(P_{w}) \cup V^{1}(P_{w})$. Since $v \in B_3$, by Corollary 1, we can assume that the $2$-colouring $b$ satisfies also the following condition:
\begin{enumerate}
\item [($a$)] if $y \in B_3$, then $b(v) \neq b(y)$.
\end{enumerate}
Let
$$M = \left\{\begin{array}{ll}
B_{3} \cup Int P_{w}, &  \mbox{if all vertices of $V^{0}(P_{w})$ are big in $G$},
\\[4pt]
B_{3}  \cup Int P_{w} \cup \{y\}, &  \mbox{if $V^{0}(P_{w}) = \{v, y\}$ and $y \in S_{3}$}.
\end{array}\right.
$$
We shall extend $b$ to a $2$-colouring $b_{0}: M  \rightarrow \{1, 2\}$ satisfying the following conditions:
\begin{enumerate}
\item [($b$)]
every cycle in $G[B_{1} \cup B_{2} \cup B_{3}] \cup G[V(P_{w}) \cup V^{0}(P_{w})]$ is not monochromatic  with respect to $a \triangledown b_{0}$,
\item [($c$)]
$b_{0}(w) = b(v)$.
\end{enumerate}
 So, we set $b_{0}(u) = b(u)$ for every $u \in B_{3}$.  It is sufficient to consider the following cases:
\begin{enumerate}
\item [($1$)]
$V^{0}(P_w) = \{v, y\}$ and $y \in B_3$,
\item [($2$)]
$V^{1}(P_w) = \{v, y\}$ and $y \in B_3$,
\item [($3$)]
$V^{1}(P_w) = \{v, y\}$ and $y \in B_1 \cup B_2$,
\item [($4$)]
$V^{0}(P_w) = \{v, y\}$ and $y \in S_3$.
\end{enumerate}

Case $(1)$. Since $y \in B_3$, by condition ($a$), we can assume that $b(v) \neq b(y)$. We  colour all vertices of $(IntP_{w})\cap~S_{1}$ with $1$ and all vertices of $(IntP_{w}) \cap S_{2}$ with $2$.  Certainly, conditions ($b$)-($c$) hold.

Case $(2)$. Since $y \in B_3$, by condition ($a$), we can assume that $b(v) \neq b(y)$.  Notice that if $w \in S_1$ ($w \in S_2$), then  $x, z \in B_2$ ($x, z \in B_1$, respectively). Thus, $b(v) \neq a(x) = a(z) = b(y)$. We colour all vertices of $IntP_{w}$ with $b(v)$.  Certainly, conditions ($b$)-($c$) hold.

Case $(3)$.
 Notice that if $w \in S_1$ ($w \in S_2$), then  $x, z \in B_2$ and $y \in B_1$ ($x, z \in B_1$ and $y \in B_2$, respectively). Thus, $b(v)= a(y) \neq a(x) = a(z)$. Moreover, there exists a vertex $s \in (IntP_{w}) \cap S_{3}$. If every  path in $G[B_{1} \cup B_{2} \cup B_{3}]$  connecting  $x$ and $z$ is not monochromatic with respect to $a \triangledown b$, then we put $b_{0}(s) = a(x)$ and we colour all  vertices of $IntP_{w} \backslash \{s\}$ with $b(v)$.

If there exists a  path in $G[B_{1} \cup B_{2} \cup B_{3}]$  connecting $x$ and $z$ which is monochromatic with respect to $a \triangledown b$, then we colour all  vertices of $IntP_{w}$ with $b(v)$.  Certainly, conditions ($b$)-($c$) hold.

Case $(4)$. We set $b_{0}(y) \neq b(v)$.  We  colour all vertices of $(IntP_{w})\cap S_{1}$ with~$1$ and all vertices of $(IntP_{w}) \cap S_{2}$ with $2$.  Certainly, conditions ($b$)-($c$) hold.

Let
$$X_{0} = B_1 \cup \{u \in M:  b_{0}(u) = 1\} \hbox{ and } Y_{0} = B_2 \cup \{u \in M:  b_{0}(u) = 2\}.$$
Remark, that $X_{0}$ and $Y_{0}$ are disjoint and, by condition (b), the induced graphs $G[X_{0}]$ and $G[Y_{0}]$ are acyclic. Hence, by Lemma  \ref{lemma1.1}, it is possible to partition the vertex set of $G$ into two subsets $S$, $T$ so that each of them induces a tree, $X_{0} \subseteq S$ and $Y_{0} \subseteq T$. Moreover,  if $w \in S_1$ ($w \in S_2$), then, by condition (c), $v, w \in X_{0}$ ($v, w \in Y_{0}$, respectively).
\end{PrfFact}

\begin{theorem}\label{theorem2.4}
Let $G \in {\cal {E}}$ and suppose that $H = G[B_{1} \cup B_{3}] \cup G[B_{2} \cup B_{3}]  \in  {\cal H}$.  If every component of $H$ is $2$-connected, then it is possible to partition the vertex set of $G$ into two disjoint subsets $S$, $T$ so that each of them induces a tree, $B_{1} \subseteq S$, $B_{2} \subseteq T$, and for every $v \in B_{3}$  the following implications are satisfied:

(1) if $d_{H}(v) \geqslant 3$, then  $N_{1}(v) \subseteq S$  and $N_{2}(v)  \subseteq T$,

(2)  if  $d_{H}(v) = 2$ and $v \in S$ $($$v \in T$$)$, then at most two vertices of  $N_{1}(v)$ $($or $N_{2}(v)$$)$ belong to  $S$  $($or $T$$)$  and all other neighbours of $v$ belong to $T$ $($$S$, respectively$)$.
\end{theorem}
 \begin{PrfFact}
Let $G \in {\cal {E}}$. Let $H = G[B_{1} \cup B_{3}] \cup G[B_{2} \cup B_{3}]  \in  {\cal H}$ and suppose that every component of $H$ is $2$-connected. Thus, $G \ncong {C^{2l}} \ast {E^2}$, for $l \geqslant 2$. Hence, every small vertex in $G$ belongs to a path of ${\cal P}_{G}$.

Notice that  $H$ is a bipartite graph with the vertex bipartition $(B_1 \cup B_2, B_3) $. Let $a: B_1 \cup B_2 \rightarrow \{1, 2\}$ be a $2$-colouring such that vertices of $B_1$ (or  $B_2$) are coloured  $1$ ($2$, respectively).  By Theorem \ref{theorem2.2}, there exists a $2$-colouring $b: B_{3} \rightarrow \{1, 2\}$ satisfying the following conditions:
\begin{enumerate}
\item[($a$)]
 every cycle in $G[B_{1} \cup B_{2} \cup B_{3}]$ is not monochromatic with respect to $a \triangledown b$,
\item[($b$)]
if a path in $H$ of length $2$ has ends belonging to $B_3$ of degree at least $3$ in $H$, and if the inner vertex is of degree $2$ in $H$, then ends of this path are coloured differently by $b$. 
\end{enumerate}
Notice that  we can also assume that $b$ satisfies the following condition:
\begin{enumerate}
\item[($c$)]
 if  a path  in $H$ of length $2$ has the inner vertex belonging to $B_3$ of degree $2$ in $H$, and if  its both ends are not coloured $2$ (are coloured $2$) by $a$, then its inner vertex is coloured $2$ ($1$, respectively) by $b$. Otherwise we can recolour the inner vertex.
\end{enumerate}

Let ${\cal R}_{G}$ be the family of all paths $P \in {\cal P}_{G}$ such that all vertices  of $V^{0}(P)$ are big in $G$ and at least one vertex of $V^{0}(P) \cup V^{1}(P)$ belongs to $B_{3}$.
Suppose that $\widehat{\cal R}_{G}$ is the family of all paths $P \in {\cal P}_{G}$ such that $V^{0}(P)$ contains a vertex belonging to $B_3$ and a vertex belonging to $S_{3}$. 

Let $P_1, \ldots, P_n$ be a sequence of all paths in ${\cal R}_{G} \cup \widehat{\cal R}_{G}$. Set $M_{0} = B_{3}$ and
$$M_{i} = \left\{\begin{array}{ll}
M_{i-1} \cup Int P_{i}, &  \mbox{for $P_{i} \in {\cal R}_{G}$},
\\[4pt]
M_{i-1} \cup Int P_{i} \cup \{y\}, &  \mbox{for $P_{i} \in {\widehat{\cal R}_{G}}$ and $y \in V^{0}(P_i) \cap S_3$}.
\end{array}\right.
$$
Denote $L_{0 } = G[B_{1} \cup B_{2} \cup B_{3}]$ and $L_{i} =  L_{i-1}  \cup  G[V(P_{i }) \cup V^{0}(P_i)]$, for $i = 1, \dots, n$.

We shall define a sequence $b_0, \ldots, b_n$ of $2$-colourings $b_{i}: M_i \rightarrow \{1, 2\}$  (we set $b_{0} = b$) such that  the following conditions, for $i > 0$, are satisfied:
\begin{enumerate}
\item[($d$)]
every cycle in $L_{i}$ is not monochromatic with respect to $a \triangledown b_i$,
\item[($e$)]
 if $v$ is a vertex belonging to $B_3 \cap (V^0(P_i) \cup V^{1}(P_i))$ with $d_{H}(v) \geqslant 3$, then all vertices
of $N_{1}(v) \cap (V(P_i) \cup V^{0}(P_i))$ (of $N_{2}(v) \cap (V(P_i) \cup V^{0}(P_i))$) are coloured $1$ ($2$, respectively) by $b_i$,
\item[($f$)]
if $v$ is a vertex belonging to $B_3 \cap (V^0(P_i) \cup V^{1}(P_i))$  with $d_{H}(v) = 2$ and if $v$ is coloured $1$ (or $2$) by $b$, then at most one vertex of $N_{1}(v) \cap  (V(P_i) \cup V^{0}(P_i))$ ($N_{2}(v) \cap (V(P_i) \cup V^{0}(P_i))$, respectively) is coloured $1$ ($2$, respectively) and all other vertices of $N(v) \cap (V(P_i) \cup V^{0}(P_i))$ are coloured $2$ ($1$, respectively) by $b_i$.
\end{enumerate}
Let $b_{i-1}$, for some $i > 0$, be a $2$-colouring satisfying conditions  ($d$)-($f$), for $i > 1$.   We shall extend $b_{i-1}$ to a $2$-colouring $b_{i}$ satisfying conditions ($d$)-($f$). So, we  put $b_{i}(u) = b_{i-1}(u)$ for every $u \in M_{i-1}$.
Let us consider a path $P_{i} \in {\cal R}_{G}$ and suppose that  $v_{i}x_{i}y_{i}z_{i}v_{i}$ is a $4$-cycle in $G$  consisting of vertices belonging to $V^{0}(P_{i}) \cup V^{1}(P_{i})$. Assume first that
$$v_{i}, y_{i} \in B_{3} \hbox{  and }x_{i}, z_{i} \in B_{1} \cup B_{2}.$$ 
Hence, $v_{i}x_{i}y_{i}z_{i}v_{i}$ is contained in $H$. Thus, it is contained in some $2$-connected component of $H$. By Lemma \ref{lem2.3}, one of the following cases occurs:
\begin{enumerate}
\item [($1$)]
$V^{0}(P_i) = \{v_{i}, y_{i}\}$, $d_{H}(v_i)$,  $d_{H}(y_i) \geqslant 3$  and $d_{H}(x_i) = deg_{H}(z_i) = 2$,
\item [($2$)]
$V^{1}(P_i) = \{v_{i}, y_{i}\}$, $d_{H}(v_i)$,  $d_{H}(y_i) \geqslant 3$ and $d_{H}(x_i) = d_{H}(z_i) = 2$,
\item [($3$)]
$V^{0}(P_i) = \{v_{i}, y_{i}\}$, $d_{H}(v_i) = d_{H}(y_i) = 2$,
\item [($4$)]
$V^{1}(P_i) = \{v_{i}, y_{i}\} $, $d_{H}(v_i) = d_{H}(y_i) = 2$.

\end{enumerate}

Case (1). Notice that $b(v_i) \neq b(y_i)$, by condition ($b$). Let us define $b_{i} (u) = 1$ ($b_{i} (u) = 2$), for every vertex $u \in S_{1} \cap Int P_i$ ($u \in S_{2} \cap Int P_i$, respectively). Certainly, conditions ($d$)-($e$) hold, for $v = v_i$ and $v = y_i$.

Case (2). Notice that $a(x_{i}) = a(z_{i})$ and $b(v_{i}) \neq b(y_{i})$, by condition ($b$). If $a(x_i) = 1$ ($a(x_i) = 2$) we set $b_{i}(u) = 2$ ($b_{i}(u) = 1$, respectively) for every vertex $u \in Int P_i$. Certainly, conditions ($d$)-($e$) hold, for $v = v_i$ and $v = y_i$.

Case (3). Notice that $b(v_i) = b(y_i)$, by  condition ($c$). Hence, $a(x_i) \neq b(v_i)$ or  $a(z_i) \neq b(v_i)$.
If  there exists a path in $L_{i-1}$ connecting $v_i$ and $y_i$ which is monochromatic with respect to $a \triangledown b_{i-1}$ and $a(x_i) \neq b(v_i)$ ($a(z_i) \neq b(v_i)$), then we  colour all inner vertices of $P_i$ with  $a(x_i)$ ($a(z_i)$, respectively).

Assume now, that every path in $L_{i-1}$ connecting $v_i$ and $y_i$  is  not monochromatic with respect to $a \triangledown b_{i-1}$. Then $a(x_i) = a(z_i) \neq b(v_i)$. If $b(v_i) = 1$  ($b(v_i) = 2$), then $x_{i}, z_{i} \in B_2$ ($x_{i}, z_{i} \in B_1$). Hence, there exists $s \in S_{1} \cap Int P_{i}$ ($s \in S_{2} \cap Int P_{i}$, respectively). We set $b_{i}(s) = 1$  ($b_{i}(s) = 2$, respectively) and we colour all other vertices of $Int P_{i}$ with $2$ ($1$, respectively). Certainly, conditions ($d$) and ($f$) hold, for $v = v_i$ and $v = y_i$.

Case (4).  Then $a(x_i) = a(z_i) \neq b(v_i) = b(y_i)$, by condition ($c$).  If there exists a path in $L_{i-1}$ connecting $x_i$ and $z_i$ which is monochromatic with respect to $a \triangledown b_{i-1}$, then we colour all vertices of $Int P_i$ with $b(v_i)$.

Assume now, that  every path in $L_{i-1}$ connecting $x_i$ and $z_i$  is not monochromatic with respect to $a \triangledown b_{i-1}$. If $b(v_i) = 1$ ($b(v_i) = 2$), then $x_{i}, z_{i}\in B_2$ ($x_{i}, z_{i} \in B_1$). Hence, there exists a vertex $s \in S_{1} \cap Int P_{i}$  ($s \in S_{2} \cap Int P_{i}$, respectively).  We set $b_{i}(s) = 2$  ($b_{i}(s) = 1$, respectively) and we colour all other vertices of $Int P_i$ with $1$ ($2$, respectively). Certainly, conditions ($d$) and ($f$) hold, for $v = v_i$ and $v = y_i$.

Assume now that
$$v_{i} \in B_{3}, y_{i} \in B_{1} \cup B_{2} \cup S_{3} \hbox{ and } x_{i}, z_{i} \in B_{1} \cup B_{2}.$$  It  is sufficient to consider the following cases:
\begin{enumerate}
\item [($5$)]
$V^{1}(P_i) = \{v_{i}, y_{i}\}$, $y_{i} \in B_1 \cup B_2$,
\item [($6$)]
$V^{0}(P_i) = \{v_{i}, y_{i}\}$, $y_{i} \in S_3$ and $d_{H}(v_i) \geqslant 3$,
\item [($7$)]
$V^{0}(P_i) = \{v_{i}, y_{i}\}$, $y_{i} \in S_3$ and $d_{H}(v_i) = 2$.
\end{enumerate}

Case $5$. Notice that $a(x_i) = a(z_i) \neq a(y_i)$.  If there exists a  path in $L_{i-1}$ connecting $x_i$ and $z_i$ which is monochromatic with respect to $a \triangledown b_{i-1}$, then we colour all vertices of $Int P_i$ with $a(y_i)$. If $d_{H}(v) \geqslant 3$, then conditions ($d$)-($e$) hold. If  $d_{H}(v) = 2$, then, by condition (c), $b(v_i) \neq a(x_i)$. Hence, $b(v_i) = a(y_i)$ and conditions ($d$) and ($f$) hold.

 Assume now,  that every path in $L_{i-1}$ connecting $x_i$ and $z_i$  is not monochromatic with respect to $a \triangledown b_{i-1}$. Thus, $b(v_i)\neq  a(x_i)$. If  $x_{i} \in B_2$ ($x_{i} \in B_1$), then  $b(v_i) = a(y_i) = 1$ ($b(v_i) = a(y_i) = 2$, respectively) and all vertices of $Int P_i$ belong to $S_{1} \cup S_{3}$ ($S_{2} \cup S_{3}$, respectively). Let $s$ be any vertex of $S_{3} \cap Int P_{i}$.  We set $b_{i}(s) = 2$ ($b_{i}(s) = 1$, respectively) and we colour all other vertices of $Int P_i$ with $1$ ($2$, respectively). Certainly, conditions ($d$)-($e$) hold.

Case (6).  If $b(v_i) = 1$ ($b(v_i) = 2$) we set $b_{i}(y_i) = 2$  ($b_{i}(y_i) = 1$, respectively). We  colour a vertex of  $Int P_{i} \cap V_1$ with $1$, and  a vertex of $Int  P_{i} \cap V_2$ with $2$. Certainly,  conditions ($d$)-($f$) hold .

Case (7).   If $b(v_i) = 1$ ($b(v_i) = 2$) we set $b_{i}(y_i) = 1$  ($b_{i}(y_i ) = 2$, respectively) and we  colour all vertices of $Int P_i$ with $2$ ($1$, respectively). Certainly,  conditions ($d$) and ($f$) hold.




Let
$$X = B_{1} \cup \{u \in M_{n}: b_{n}(u) = 1\} \hbox{ and } Y = B_{2} \cup \{u \in M_{n}: b_{n}(u) = 2\}.$$
Remark, that $X$ and $Y$ are disjoint and, by condition (d),  graphs $G[X]$ and $G[Y]$ are acyclic. For every path in ${\cal R}_{G}\cup \widehat{\cal R}_{G}$ (in ${\cal P}_G \backslash ({\cal R}_{G}\cup \widehat{\cal R}_{G})$) the set of its inner vertices is contained in $X \cup Y$ (is disjoint with  $X \cup Y$). Therfore, by  Lemma  \ref{lemma1.1}, it is possible to partition the vertex set of $G$ into two subsets $S$, $T$ so that each of them induces a tree, $X \subseteq S$ and $Y \subseteq T$.

Notice that every facial cycle of the graph $G[B_{1} \cup B_{2} \cup B_{3} ]$ is a square or a triangle. Moreover, if  $d_{H}(v) = 2$, then there are exactly two facial cycles containing the vertex $v$.  Hence, by conditions ($e$)-$(f$) (for every $i = 1, \dots, n$) the following implications hold:
\begin{enumerate}
\item[($g$)]
if $v \in B_3$ has degree at least $3$ in $H$, then  vertices of $N_{1}(v)$ (of $N_{2}(v)$) are coloured $1$ ($2$, respectively)  by $b_n$,
\item[$(h$)]
 if  $v \in B_3$ has degree $2$ in  $H$ and it is coloured $1$ (or $2$) by $b$, then at most two vertices of  $N_{1}(v)$ (of $N_{2}(v)$) are coloured $1$ (or $2$),  and all other neighbours of $v$ are coloured $2$ ($1$, respectively) by $b_n$.
\end{enumerate}
Finally, by conditions ($g$)-$(h$),  Theorem \ref{theorem2.4} holds.
\end{PrfFact}

\end{document}